\documentstyle[leqno,11pt]{article}
\def\li{\mathop{\mathrm{li}}}

\def\eps{\varepsilon}

\begin{document}

\title{\sc On certain arithmetic functions involving exponential divisors}
\author{{\bf L\'aszl\'o T\'oth} (P\'ecs, Hungary)}
\date{ }
\maketitle

\centerline{Annales Univ. Sci. Budapest., Sect. Comp., {\bf 24}
(2004), 285-294} \vskip4mm

\abstract{The integer $d$ is called an exponential divisor of
$n=\prod_{i=1}^r p_i^{a_i}>1$ if $d=\prod_{i=1}^r p_i^{c_i}$, where
$c_i | a_i$ for every $1\le i \le r$. The integers $n=\prod_{i=1}^r
p_i^{a_i}, m=\prod_{i=1}^r p_i^{b_i}>1$ having the same prime
factors are called exponentially coprime if $(a_i,b_i)=1$ for every
$1\le i\le r$.

In this paper we investigate asymptotic properties of certain
arithmetic functions involving exponential divisors and
exponentially coprime integers.}

\vskip5mm

{\bf 1. Introduction} \vskip3mm

Let $n>1$ be an integer of canonical form $n=\prod_{i=1}^r p_i^{a_i}$.
The integer $d$ is called an {\sl exponential divisor} of $n$ if
$d=\prod_{i=1}^r p_i^{c_i}$, where $c_i | a_i$ for every $1\le i \le r$, notation: $d|_e n$.
By convention $1|_e 1$. This notion was introduced by {\sc M. V. Subbarao} \cite{Su72}.
Note that $1$ is not an exponential divisor of $n>1$, the smallest exponential divisor of $n>1$ is its squarefree kernel $\kappa(n)=\prod_{i=1}^r p_i$.

Let $\tau^{(e)}(n)= \sum_{d|_e n} 1$ and $\sigma^{(e)}(n)=\sum_{d|_e n} d$ denote the number and the sum of exponential divisors of $n$, respectively.
The integer $n=\prod_{i=1}^r p_i^{a_i}$ is called {\sl exponentially squarefree} if all the exponents $a_i$ ($1\le i \le r$) are squarefree. Let $q^{(e)}$ denote the characteristic function of exponentially squarefree integers. Properties of these functions were investigated by several authors, see \cite{FaSu89}, \cite{KaSu99}, \cite{KaSu2003}, \cite{PeWu97}, \cite{SmWu97}, \cite{Su72}, \cite{Wu95}.

Two integers $n,m >1$ have common exponential divisors iff they have the same prime
factors and in this case, i.e. for $n=\prod_{i=1}^r p_i^{a_i}$, $m=\prod_{i=1}^r p_i^{b_i}$, $a_i,b_i\ge 1$ ($1\le i\le r$), the {\sl greatest common exponential divisor} of $n$ and $m$ is
$$
(n,m)_{(e)}:=\prod_{i=1}^r p_i^{(a_i,b_i)}.
$$

Here $(1,1)_{(e)}=1$ by convention and $(1,m)_{(e)}$ does not exist for $m>1$.

The integers $n,m >1$ are called {\sl exponentially coprime}, if
they have the same prime factors and $(a_i,b_i)=1$ for every $1\le
i\le r$, with the notation of above. In this case $(n,m)_{(e)}=\kappa(n)=\kappa(m)$.
$1$ and $1$ are considered to be exponentially coprime. $1$ and
$m>1$ are not exponentially coprime.

For $n=\prod_{i=1}^r p_i^{a_i}>1$, $a_i\ge 1$ ($1\le i\le r$), denote by
$\phi^{(e)}(n)$ the number of integers $\prod_{i=1}^r p_i^{c_i}$ such that $1\le c_i\le a_i$ and $(c_i,a_i)=1$ for $1\le i \le r$, and let $\phi^{(e)}(1)=1$.
Thus $\phi^{(e)}(n)$ counts the number of divisors $d$ of $n$ such that $d$ and $n$ are
exponentially coprime.

It is immediately, that $\phi^{(e)}$ is a prime independent multiplicative function and for $n>1$,
$$
\phi^{(e)}(n)=\prod_{i=1}^r \phi(a_i),
$$
where $\phi$ is the Euler-function.
Exponentially coprime integers and function $\phi^{(e)}$
were introduced by {\sc J. S\'andor} \cite{Sa96}. He showed that
$$
\limsup_{n\to \infty} \frac{\log \phi^{(e)}(n) \log \log n}{\log n}=
\frac{\log 4}{5}.
$$

We consider the functions $\tilde{\sigma}$ and $\tilde{P}$ defined as follows. Let $\tilde{\sigma}(n)$ be the sum of those divisors $d$ of $n$ such that $d$ and $n$ are exponentially coprime. Function $\tilde{\sigma}$ is multiplicative and for every prime power $p^a$,
$$
\tilde{\sigma}(p^a) =\sum_{1\le c\le a \atop (c,a)=1} p^c.
$$
Here $\tilde{\sigma}(p)=\tilde{\sigma}(p^2)=p, \tilde{\sigma}(p^3)=p+p^2, \tilde{\sigma}(p^4)=p+p^3$, etc.

Furthermore let $\tilde{P}(n)$ be given by
$$
\tilde{P}(n)=\sum_{1\le j\le n \atop \kappa(j)=\kappa(n)} (j,n)_{(e)},
$$
representing an analogue of Pillai's function $P(n)=\sum_{j=1}^n (j,n)$.

Function $\tilde{P}$ is also multiplicative and for every prime power $p^a$,
$$
\tilde{P}(p^a) = \sum_{1\le c\le a} p^{(c,a)}= \sum_{d | a} p^d \phi(a/d),
$$
here $\tilde{P}(p)=p, \tilde{P}(p^2)=p+p^2, \tilde{P}(p^3)=2p+p^3, \tilde{P}(p^4)=2p+p^2+p^4$, etc.

We call an integer $n=\prod_{i=1}^r p_i^{a_i}$ {\sl exponentially $k$-free} if all the exponents $a_i$ ($1\le i \le r$) are $k$-free, i.e. are not divisible by the $k$-th power of any prime ($k\ge 2$). Let $q_k^{(e)}$ denote the characteristic function of exponentially $k$-free integers.

The aim of this paper is to investigate the functions $\phi^{(e)}(n)$, $\tilde{\sigma}(n)$,
$\tilde{P}(n)$ and $q_k^{(e)}(n)$. The estimate given for the sum $\sum_{n\le x} q_k^{(e)}(n)$ generalizes the result of {\sc J. Wu} \cite{Wu95} concerning exponentially squarefree integers. Our main results are formulated in Section 2, their proofs are given in Section 3.

Our estimates for $\sum_{n\le x} (\tilde{\sigma}(n))^u$ and $\sum_{n\le x} q_k^{(e)}(n)$ are consequences of a general result due to {\sc V. Sita Ramaiah} and {\sc D. Suryanarayana} \cite{SiSu82}, the proof of which uses the estimate  of {\sc A. Walfisz} \cite{Wa63} concerning $k$-free integers and is simpler than the proof given by {\sc J. Wu} \cite{Wu95}.

{\sc A. Smati} and {\sc J. Wu} \cite{SmWu97} deduced some interesting analogues of
known results on the divisor function $\tau(n)$ in case of $\tau^{(e)}(n)$. They remarked
that their results can be stated also for certain other prime independent multiplicative functions $f$ if $f(n)$ depends only on the squarefull kernel of $n$.

We point out two such results in case of $\phi^{(e)}(n)$. Note that, since $\phi(1)=\phi(2)=1$, $\phi^{(e)}(n)$ depends only on the cubfull kernel of $n$.
These results are contained in Section 4. Here some open problems are also stated.

\vskip3mm
{\bf 2. Main results}
\vskip3mm

Regarding the average orders of the functions $\phi^{(e)}(n)$, $\tilde{\sigma}(n)$
and $\tilde{P}(n)$ we prove the following results.

\vskip2mm
{\bf Theorem 1.} {\it
$$
\sum_{n\le x} \phi^{(e)}(n)=C_1x +C_2 x^{1/3} + O(x^{1/5+\eps}),
$$
for every $\eps >0$, where $C_1, C_2$ are constants given by
$$
C_1=\prod_p \left(1 + \sum_{a=3}^{\infty} \frac{\phi(a)-\phi(a-1)}{p^a}\right),
$$
$$
C_2=\zeta(1/3) \prod_p \left(1 + \sum_{a=5}^{\infty} \frac{\phi(a)-\phi(a-1)-\phi(a-3)+
\phi(a-4)}{p^{a/3}}\right).
$$}
\vskip2mm

\vskip2mm
{\bf Theorem 2.} {\it Let $u >1/3$ be a fixed real number. Then
$$
\sum_{n\le x} (\tilde{\sigma}(n))^u = C_3 x^{u+1} +O(x^{u+1/2}\delta(x)),
$$
where $C_3$ is given by
$$
C_3=\frac1{u+1} \prod_p \left(1+ \sum_{a=2}^{\infty} \frac{(\tilde{\sigma}(p^a))^u-p^u (\tilde{\sigma}(p^{a-1}))^u}{p^{a(u+1)}}\right)
$$
and
$$
\delta(x)=\exp(-A(\log x)^{3/5} (\log \log x)^{-1/5}),
$$
$A$ being a positive constant. }
\vskip2mm

\vskip2mm
{\bf Theorem 3.} {\it
$$
\sum_{n\le x} \tilde{P}(n) = C_4 x^2 + O(x (\log x)^{5/3}),
$$
where the constant $C_4$ is given by
$$
C_4=\frac1{2} \prod_p \left(1+ \sum_{a=2}^{\infty} \frac{\tilde{P}(p^a)-p \tilde{P}(p^{a-1})}{p^{2a}}\right).
$$ }
\vskip2mm

Concerning the maximal order of the function $\tilde{P}(n)$ we have

\vskip2mm
{\bf Theorem 4.} {\it
$$
\limsup_{n\to \infty} \frac{\tilde{P}(n)}{n\log\log n}
= \frac{6}{\pi^2}e^{\gamma},
$$
where $\gamma$ is Euler's constant.}
\vskip2mm

\vskip2mm
{\bf Theorem 5.} {\it If $k\ge 2$ is a fixed integer, then
$$
\sum_{n\le x} q_k^{(e)}(n) = D_k x +O(x^{1/2^k} \delta(x)),
$$
where
$$
D_k=\prod_p \left(1+ \sum_{a=2^k}^{\infty} \frac{q_k(a)-q_k(a-1)}{p^a}\right),
$$
$q_k(n)$ denoting the characteristic function of $k$-free integers.}
\vskip2mm

In the special case $k=2$ case this formula is due to {\sc J. Wu} \cite{Wu95}, improving an earlier result of {\sc M. V. Subbarao} \cite{Su72}.
\newpage

{\bf 3. Proofs}
\vskip3mm

The proof of Theorem 1 is based on the following lemma.

\vskip2mm
{\bf Lemma 1.} {\it The Dirichlet series of $\phi^{(e)}$ is absolutely convergent for
$Re\; s >1$ and it is of form
$$
\sum_{n=1}^{\infty} \frac{\phi^{(e)}(n)}{n^s} =\zeta(s)\zeta(3s) V(s),
$$
where the Dirichlet series $V(s)=\sum_{n=1}^{\infty} \frac{v(n)}{n^s}$ is absolutely convergent for $Re\; s > 1/5$. }
\vskip2mm

\vskip2mm
{\bf Proof of Lemma 1.} Let $\mu_3(n)=\mu(m)$ or $0$, according as
$n=m^3$ or not, where $\mu$ is the M\"obius function, and let $f=\mu_3*\mu$ in terms of the Dirichlet convolution. Then we can formally obtain the desired expression by taking $v=\phi^{(e)}*f$.  Both $f$ and $v$ are multiplicative and easy computations show that
$f(p)=f(p^3)=-1, f(p^4)=1, f(p^2)=f(p^a)=0$ for each $a\ge 5$, and $v(p^a)=0$ for $1\le a \le 4$, $v(p^a)=\phi(a)-\phi(a-1)-\phi(a-3)+\phi(a-4)$ for $a\ge 5$.

Since $|v(p^a)|< 4a$ for $a\ge 5$, we obtain that $V(s)$ is absolutely convergent for $Re\; s > 1/5$.
\vskip2mm

\vskip2mm
{\bf Proof of Theorem 1.} Lemma 1 shows that $\phi^{(e)}=v*\tau(1,3,\cdot)$, where
$\tau(1,3,n)=\sum_{ab^3=n} 1$ for which
$$
\sum_{n\le x} \tau(1,3,n) =\zeta(3)x+\zeta(1/3)x^{1/3}+O(x^{1/5}),
$$
cf. \cite{Kr88}, p. 196-199. Therefore,
$$
\sum_{n\le x} \phi^{(e)}(n) =\sum_{d\le x} v(d) \sum_{e \le x/d} \tau(1,3,e) =
$$
$$
=\zeta(3)x \sum_{d\le x} \frac{v(d)}{d}+ \zeta(1/3)x^{1/3}
\sum_{d\le x} \frac{v(d)}{d^{1/3}}+O\left(x^{1/5+\eps}\sum_{d\le x}
\frac{|v(d)|}{d^{1/5+\eps}} \right),
$$
and obtain the desired result by usual estimates.
\vskip2mm

For the proof of Theorem 2 we use the following general result due to {\sc V. Sita Ramaiah} and {\sc D. Suryanarayana} \cite{SiSu82}, Theorem 1.

\vskip2mm
{\bf Lemma 2.} {\it Let $k\ge 2$ be a fixed integer, $\beta >(k+1)^{-1}$ be a fixed real number and $g$ be a multiplicative arithmetic function such that $|g(n)|\le 1$ for all $n\ge 1$. Suppose that either

(i) $|g(p^j)-1|\le p^{-1}$ for $1\le j \le k-1$, $g(p^k)=0$ for all primes $p$, or

(ii) $g(p^j)=1$ for $1\le j \le k-1$, $g(p^k)=p^{-\beta}$ for all primes $p$.

Then
$$
\sum_{n\le x} g(n)= x \sum_{n=1}^{\infty} \frac{(g*\mu)(n)}{n}+ O(x^{1/k}\delta(x)).
$$ }
\vskip2mm

\vskip2mm
{\bf Proof of Theorem 2.} This is a direct consequence of Lemma 2 of above.
Take $g(n)=(\tilde{\sigma}(n)/n)^u$. Here $g(p)=1$, $g(p^2)=p^{-u}$, $g(p^a)\le
p^{-au}(p+p^2+...+p^{a-1})^u < (p-1)^{-u} \le 1$ for every $a\ge 3$, hence $0<g(n)\le 1$ for all $n\ge 1$. Choosing $k=2$, $\beta =u$, we obtain the given result by partial summation.
\vskip2mm

\vskip2mm
{\bf Lemma 3.} The Dirichlet series of $\tilde{P}(n)$ is absolutely convergent for
$Re\; s >2$ and it is of form
$$
\sum_{n=1}^{\infty} \frac{\tilde{P}(n)}{n^s} =\frac{\zeta(s-1)\zeta(2s-1)}{\zeta(3s-2)} W(s),
$$
where the Dirichlet series $W(s)=\sum_{n=1}^{\infty} \frac{w(n)}{n^s}$ is absolutely convergent for $Re\; s > 3/4$.
\vskip2mm

\vskip2mm
{\bf Proof of Lemma 3.}
$$
\sum_{n=1}^{\infty} \frac{\tilde{P}(n)}{n^s} =\prod_p \left(1+\sum_{a=1}^{\infty} \sum_{d|a} \frac{p^d\phi(a/d)}{p^{as}}\right)
$$ $$
=\prod_p \left(1+\sum_{j=1}^{\infty} \phi(j) \sum_{d=1}^{\infty} \frac1{p^{d(js-1)}}\right)=
\prod_p \left(1+\sum_{j=1}^{\infty} \frac{\phi(j)}{p^{js-1}-1}\right)
$$
$$
=\frac{\zeta(s-1)\zeta(2s-1)}{\zeta(3s-2)}W(s),
$$
where
$$
W(s):= \prod_p \left(1 + \frac{(p^{s-1}-1)
(p^{2s-1}-1)}{p^{3s-2}-1} \sum_{j=3}^{\infty} \frac{\phi(j)}{p^{js-1}-1}\right),
$$
which is absolutely convergent for $Re\; s > 3/4$.
\vskip2mm

\vskip2mm
{\bf Proof of Theorem 3.} By Lemma 3, $\tilde{P}= h*w$, where
$$
h(n)=\sum_{ab^2c^3=n} abc^2\mu(c),
$$
and obtain the desired result, exactly like in proof of Theorem 2 of \cite{PeWu97},
using the estimate
$$
\sum_{mn^2\le x} mn=\frac1{2}\zeta(3)x^2 +O(x(\log x)^{2/3})
$$
due to {\sc Y. - F. S. P\'etermann} and {\sc J. Wu} \cite{PeWu97}, Theorem 1.
\vskip2mm

Theorem 4 is a direct consequence of the following general result of {\sc L. T\'oth} and {\sc E. Wirsing} \cite{ToWi2003}, Corollary 1.

\vskip2mm
{\bf Lemma 4.} {\it Let $f$ be a nonnegative real-valued multiplicative function. Suppose that
for all primes $p$ we have $\rho(p):=\sup_{\nu \ge 0}f(p^\nu) \le (1-1/p)^{-1}$ and that for all primes $p$ there is an exponent $e_p=p^{o(1)}$ such that $f(p^{e_p})\ge 1+1/p$. Then
$$
\limsup_{n\to \infty} \frac{f(n)}{\log \log n}=e^{\gamma}\prod_p \left(1-\frac1{p}\right) \rho(p).
$$}
\vskip2mm

\vskip2mm
{\bf Proof of Theorem 4.} Apply Lemma 4 for $f(n)=\tilde{P}(n)/n$, where $f(p^a)\le (p+p^2+\cdots +p^a)p^{-a}<(1-1/p)^{-1}$ for every $a\ge 1$ and $f(p^2)=1+1/p$, hence we can choose $e_p=2$ for all $p$. Moreover, $\rho(p)=1+1/p$ for all $p$ and obtain the desired result.
\vskip2mm

\vskip2mm
{\bf Proof of Theorem 5.}
This follows from Lemma 2 by taking $2^k$ instead of $k$, where $q_k^{(e)}(p)=q_k^{(e)}(p^2)=...=q_k^{(e)}(p^{2^k-1})=1$, $q_k^{(e)}(p^{2^k})=0$.
\vskip2mm

\vskip3mm
{\bf 4. Further results and problems}
\vskip3mm

The next result is an analogue of the exponential divisor problem of Titchmarsh, see
Theorem 1 of \cite{SmWu97}. The proof is the same using that $\phi^{(e)}(n)$ is a prime independent multiplicative function depending only on the squarefull (cubfull) kernel of $n$
and that $\phi^{(e)}(p^a)=\phi(a)\le a$ for every $a\ge 1$.

\vskip2mm
{\bf Theorem 6.} {\it For every fixed $B>0$,
$$
\sum_{p\le x} \phi^{(e)}(p-1)=C_5 \li x + O(x/(\log x)^B),
$$
where
$$
C_5=\prod_p \left(1 + \sum_{k=3}^{\infty} \frac{\phi(k)-1}{p^k} \right).
$$}
\vskip2mm

Let $\omega(n)$ and $\Omega(n)$ denote, as usual, the number of prime factors of $n$ and the number of prime power factors of $n$, respectively.

\vskip2mm
{\bf Theorem 7.} {\it A maximal order of  $\Omega(\phi^{(e)}(n))$ is $2(\log n)/5\log \log n$. }
\vskip2mm

This can be obtained by the same arguments as those given in the proof of Theorem 3.(i) of \cite{SmWu97}. Here the upper bound is attained for $n_k=(p_1\cdots p_k)^5$, where $p_k$ is the $k$-th prime.

\vskip2mm
{\bf Problem 1.} {\it Determine a maximal order of $\omega(\phi^{(e)}(n))$.}
\vskip2mm

Since $\tilde{\sigma}(n) \le n$ for all $n\ge 1$ and $\tilde{\sigma}(p)=p$ for all primes $p$, it is clear that a maximal order of $\tilde{\sigma}(n)$ is $n$.

\vskip2mm
{\bf Problem 2.} {\it Determine a minimal order of $\tilde{\sigma}(n)$.}
\vskip2mm

{\sc J. S\'andor} \cite{Sa96} considered in fact the function $\varphi_e(n)$ defined
as the number of integers $1<a<n$ for which $a$ and $n$ are exponentially coprime ($n>1$) and $\varphi_e(1)=1$.
Although $\varphi_e(p^a)=\phi^{(e)}(p^a)=\phi(a)$ for any prime power $p^a$, functions $\varphi_e$ and $\phi^{(e)}$ are not the same. Take for example $n=2^3\cdot 3^2$, then numbers $a<n$ exponentially coprime to $n$ are $a=2\cdot 3, 2^2\cdot 3, 2^4\cdot 3$,
hence $\varphi_e(2^3 \cdot 3^2) = 3 \ne 2 \cdot 1 =\phi(3)\phi(2)=\varphi_e(2^3) \cdot \varphi_e(3^2)$.

Therefore, $\varphi_e$ is not multiplicative and $\varphi_e(n)\ge \phi^{(e)}(n)$ for every $n\ge 1$.

\vskip2mm
{\bf Problem 3.} {\it What can be said on the order of the function $\varphi_e(n)$?}

\vskip4mm

\noindent{{\bf L\'aszl\'o T\'oth}\\
University of P\'ecs\\
Institute of Mathematics and Informatics\\
Ifj\'us\'ag u. 6\\
7624 P\'ecs, Hungary\\
ltoth@ttk.pte.hu}

\end{document}